\documentclass[11pt]{amsart}

\issueinfo{00}{0}{Month}{Year}
\copyrightinfo{2008}{University of Calgary}
\pagespan{1}{2}
\PII{ISSN 1715-0868}
\usepackage{graphicx}
\usepackage{epstopdf}
\usepackage{amsthm,amsmath,amssymb,amsxtra,amscd}
\usepackage{verbatim}
\usepackage{indent}
\usepackage{rotating}
\usepackage{multirow}
\usepackage{multicol}
\usepackage{url}
\usepackage{algorithmic}
\usepackage{xypic}

\newtheorem{theorem}{Theorem}[section]

\newtheorem{lemma}[theorem]{Lemma}

\newtheorem{definition}[theorem]{Definition}

\newtheorem*{example*}{Example}

\newtheoremstyle{myexample}{3pt}{3pt}{\rmfamily}{}{\itshape}{:}{ }{\thmname{#1}\thmnumber{ #2}\thmnote{ (#3)}}
\theoremstyle{myexample}

\newtheoremstyle{myremark}{3pt}{3pt}{\rmfamily}{}{\itshape}{:}{ }{\thmname{#1}}
\theoremstyle{myremark}

\newtheorem*{observation*}{Observation}

\newtheoremstyle{conjecture}{3pt}{3pt}{\itshape}{}{\bfseries}{.}{ }{\thmname{#1}\thmnote{ (#3)}}
\theoremstyle{conjecture}

\newtheorem*{question*}{Question}

\newtheorem{theorem*}{Theorem}

\numberwithin{equation}{section}

\newcounter{algorithm}
\renewcommand{\thealgorithm}{\thesection.\arabic{algorithm}}

\usepackage{graphicx}
\usepackage{subcaption}

\begin{document}

% A short title is not required, but if needed use:
% \title[short title]{full title}
\title{Two properties of maximal antichains in strict chain product posets}

\author{Shen-Fu Tsai}
\address{Google LLC, 747 6th Street South, Kirkland, WA, USA}
\email{parity@gmail.com; parity@google.com}

% For each additional author, add another set of
% \author, \address, and \email commands
%
%\author{}
%\address{}
%\email{}

\date{(date1), and in revised form (date2).}

% \thanks entries are to acknowledge grants. You may combine
% all acknowledgments into one \thanks entry, or may use
% multiple \thanks entries. They generate footnotes without
% tags, so you must be explicit about which authors are
% thanking whom.
\thanks{The author would like to thank the anonymous reviewers for their valuable comments
and suggestions, and Jesse Geneson for improving the note's readability.}
\thanks{This work is not related to the author's activities at Google.}

% Insert paper text here
\begin{abstract}
We present two results on maximal antichains in the strict chain product poset $[t_1+1]\times[t_2+1]\times\ldots\times[t_n+1]$. First, we prove that these maximal antichains are also maximum. Second, we prove that there is a bijection between maximal antichains in the strict chain product poset $[t_1+1]\times[t_2+1]\times\ldots\times[t_n+1]$ and antichains in the non-strict chain product poset $[t_1]\times[t_2]\times\ldots\times[t_n]$.
\end{abstract}

\maketitle

\section{Introduction}
\label{intro}
A \textit{partially ordered set} (poset) is a set $S$ together with a reflexive, antisymmetric transitive relation, often denoted by $\leq$. In a \textit{strict} or \textit{irreflexive} poset, the relation is irreflexive and transitive and is denoted by $<$. Call $S$ \textit{totally ordered} provided that for all distinct $x, y$ in $S$, $x < y$ or $y < x$. A \textit{direct product} of posets $S$ and $T$ is the partially ordered set defined on the cartesian product $S \times T$ with its relation defined componentwise by the orders of $S$ and $T$. When the factors of a direct product are chains, call the resulting poset a \textit{chain product poset}. An \textit{antichain} in $S$ has no distinct elements that are comparable.

In this note we show that all antichains that are maximal (with respect to set containment) in any strict chain product poset are also maximum (that is, of maximum cardinality among all antichains). We also prove that there is a bijection between maximal antichains in any strict chain product poset and all antichains in a naturally-related (non-strict) chain product poset.

\begin{definition}
Given a positive integer $t$, let $[t]$ denote $\{1,2,\ldots,t\}$. For tuples $x=\left(x_1,\ldots,x_n\right)$ and $y=\left(y_1,\ldots,y_n\right)$, we say that $x\leq y$ if $x_i\leq y_i$ for every $i\in[n]$, and $x<y$ if $x_i< y_i$ for every $i\in[n]$. The chain product poset $[t_1]\times\cdots\times[t_n]$ is a poset containing $\prod_{i\in[n]}t_i$ tuples, each represented by a $n$-tuple $\left(x_1,\ldots,x_n\right)$ where $x_i\in[t_i]$ for every $i\in[n]$.\\
\end{definition}

\section{Results}
\begin{lemma}
\label{lemma:maximal_equivalent_to_maximum}
For any strict poset $\mathcal{G}=[t_1+1]\times[t_2+1]\times\cdots\times[t_n+1]$, the maximal antichains in $\mathcal{G}$ are maximum, and each has size $\prod_{i\in[n]}\left(t_i+1\right)-\prod_{i\in[n]}t_i$~\cite{tsai_maa_2017}.
\end{lemma}
We first present a useful notion: shape. \footnote{Our original proof was long, yet the reviewers of our previous attempted submission of Lemma~\ref{lemma:maximal_equivalent_to_maximum} proposed the concept of {\it shape} and shortened its proof. We cannot be more grateful for their help.\\}
\begin{definition}
	A {\it shape} is a tuple $s=\left(s_1,\ldots,s_n\right)$ such that $\min_{i\in[n]}s_i=1$. A tuple $x=\left(x_1,\ldots,x_n\right)$ has shape $s$ if for some $h$, $x_i=s_i+h$ for every $i\in[n]$. Given a set of tuples $T$ and a shape $s$, denote the set of tuples in $T$ with shape $s$ by $P\left(T,s\right)$.\\
\end{definition}

\begin{proof}
There are $\prod_{i\in[n]}\left(t_i+1\right)-\prod_{i\in[n]}t_i$ shapes contained in $\mathcal{G}$, call this set of shapes $S$. For every $s\in S$, each antichain in $\mathcal{G}$ contains at most one tuple of shape $s$, so it suffices to show $s$ is present in any maximal antichain in $\mathcal{G}$.\\

	Suppose that $s$ has no presence in some maximal antichain $\mathcal{M}\subset\mathcal{G}$. Let the tuples in $\mathcal{G}$ with shape $s$ be $z_1<z_2<\ldots<z_k$ where for every $i\in[k-1]$, $z_{i+1}-z_i=1$. If $k=1$, then no other tuple in $\mathcal{G}$ is comparable with $z_1$ and therefore $z_1\in \mathcal{M}$.  Otherwise there exists $i\in[k-1]$ such that $z_i<x\in \mathcal{M}$ and $z_{i+1}>y\in \mathcal{M}$. Notice that $y\leq z_i<z_{i+1}\leq x$, a contradiction.\\
\end{proof}

\begin{lemma}
	There is a bijection between antichains in the non-strict poset $\mathcal{F}=[t_1]\times[t_2]\times\cdots\times[t_n]$ and maximal antichains in the strict poset $\mathcal{G}=[t_1+1]\times[t_2+1]\times\cdots\times[t_n+1]$.\\
\end{lemma}

Given a non-strict poset $\mathcal{P}$, an \textit{(order) ideal} in $\mathcal{P}$ is a subset $\mathcal{I}\subset \mathcal{P}$ such that if $x\in \mathcal{I}$, $y\in \mathcal{P}$, and $y\leq x$, then $y\in \mathcal{I}$. An ideal is uniquely characterized by its maximal elements, so to prove the lemma it suffices to find a bijection between ideals in $\mathcal{F}$ and maximal antichains in $\mathcal{G}$.\

\begin{proof}
The mapping $\phi$ below is the bijection that serves our purpose.\\

Define the function $H\left(x\right)=x+1$ for any tuple $x$. Given an ideal $\mathcal{I}\subset\mathcal{F}$, define $\phi\left(\mathcal{I}\right)=\phi_1\left(\mathcal{I}\right)\cup\phi_2\left(\mathcal{I}\right)$. For every shape $s$ present in $H\left(\mathcal{I}\right)$, include in $\phi_1\left(\mathcal{I}\right)$ the maximal tuple from $H\left(\mathcal{I}\right)$ that has shape $s$. For every other shape present in $\mathcal{G}$, include in $\phi_2\left(\mathcal{I}\right)$ the smallest tuple from $\mathcal{G}$ of that shape. Sets $\phi_1\left(\mathcal{I}\right)$ and $\phi_2\left(\mathcal{I}\right)$ are disjoint.\\

{\bf Injectivity}. Given two distinct ideals $\mathcal{I}_1$ and $\mathcal{I}_2$, if there exists a shape $s$ present in $H\left(\mathcal{I}_1\right)$ but not $H\left(\mathcal{I}_2\right)$, then
$$\min\left(\min\left(P\left(\phi\left(\mathcal{I}_2\right),s\right)\right)\right)=1<2\leq\min\left(\min\left(P\left(\phi\left(\mathcal{I}_1\right), s\right)\right)\right).$$
Otherwise there exists a shape $s$ such that $P\left(H\left(\mathcal{I}_1\right), s\right)\neq P\left(H\left(\mathcal{I}_2\right), s\right)$ and both are not empty. Since they are both ideals in $\left([t_1+1]-\{1\}\right)\times\cdots\times\left([t_n+1]-\{1\}\right)$,
	$$
	P\left(\phi\left(\mathcal{I}_1\right), s\right)=\{\max\left(P\left(H\left(\mathcal{I}_1\right), s\right)\right)\}\neq \{\max\left(P\left(H\left(\mathcal{I}_2\right), s\right)\right)\}=P\left(\phi\left(\mathcal{I}_2\right), s\right).
	$$
\\
{\bf Surjectivity}.
	Every maximal antichain $\mathcal{M}$ in $\mathcal{G}$ is the image of some ideal in $\mathcal{F}$ under $\phi$: let the maximal tuples of an ideal $\mathcal{I}$ in $\mathcal{F}$ be $\{x-1|x\in\mathcal{M}, x-1\geq 0\}$. Clearly $\phi(\mathcal{I})=\mathcal{M}$.\\

It remains to prove that $\phi\left(\mathcal{I}\right)$ is an antichain in $\mathcal{\mathcal{G}}$. By Lemma~\ref{lemma:maximal_equivalent_to_maximum}, this implies its  maximality as each shape present in $\mathcal{G}$ is also present in $\phi\left(\mathcal{I}\right)$.\\

Suppose in $\phi\left(\mathcal{I}\right)$ that there exist tuples $x<y$. Tuple $y$ belongs to $\phi_1\left(\mathcal{I}\right)$ because every tuple in $\phi_2\left(\mathcal{I}\right)$ has an element $1$. Tuple $x\notin\mathcal{I}$ because otherwise $x+k\in\phi\left(\mathcal{I}\right)$ for some positive integer $k$. However $x=H^{-1}\left(x+1\right)\leq H^{-1}\left(y\right)\in\mathcal{I}$ contradicts $\mathcal{I}$ being an ideal in $\mathcal{F}$.\\
\end{proof}

\nocite{*}
\bibliographystyle{amsplain}
\bibliography{cdm_2020_04}
\end{document}